\documentclass[11pt]{article}
\usepackage{amsmath}
\usepackage{amssymb}
\begin{document}
\title{An  Index  Interpretation For  The  Number Of  Limit  Cycles  of a  Vector  Field}
\author{Ali Taghavi \\ Institute for Advanced Studies in Basic Sciences \\  Zanjan 45195-159, Iran}
\maketitle

Does the Hilbert 16th Problem  have   a  $PDE$  nature?The Hilbert
16th
 Problem asks for a uniform  upper bound  $H(n)$ for the  number  of  limit  cycles of a polynomial  vector
 field $X$  of  degree  $n$  on  the  Plane. More ever   it  seems that  "limit  cycles"
 are  The  only  obstructions  for  solving  the  $"PDE"$ $X.g=f$ ,  globally  in the plane

 The following  observation about  Lienard  equation suggests  to  look at the Hilbert 16th problem
 as  a  $PDE$  problem\\\\
\textbf{ Proposition.} Let L be the Lienard polynomial Vector
field
  $$
\begin{cases}
 \dot{x}=y-F(x)\\
 \dot{y}=-x
\end{cases}
$$
where  F is  an odd degree  polynomial with $F'(0)\neq0$ ,L
defines  a  linear operator on  function space, $L(f)=L.f$  If All
Limit  Cycles  of L be Hyperbolic  Then The  number of  limit
cycles of  $L$ is equal  to  $-index L$.\\
\textbf{Proof.}The origin is The  only  singularity  of  vector
field $L$ and  let we have n limit
cycles,$\gamma_{1}$,$\gamma_{2}$\ldots$\gamma_{n}$ which all
surround the origin  Let  $f$  satisfies  the following
conditions: its integral along all closed orbits is zero and
$f(0)=0$, actually such f is in the kernel of an operator defined
on the  function space
to $n+1$ dimension \\
We prove there is a map g with $L.g=f$,This shows that the
Fredholm Index of L is equal to $-n$,because the kernel is one
dimensional space since around attractors the only first integrals
are constant maps ).\\ Since the origin is a Hyperbolic
singularity,we can define g in a unique integral way in the
interior
 of $\gamma_{1}$, see below as a similar situation near hyperbolic  limit  cycle ,g is  uniform continuous  in the interior of
 $\gamma_{1}$ and has a  unique  extension to the boundary, because the
 integral of f along $\gamma_{1}$ is  zero and $g(x)-g(p(x)$ is
 near to zero  where  $p(x)$ is a poincare  map with respect to
 some transverse section.Now We extend  g to exterior of
 $\gamma_{1}$ s  follows:
 We  Define $g(x)=g(x^{*})+\int_{0}^\infty
 f(\varphi_{t}(x))-f(\varphi_{t}(x^{*}))$.
 $x^{*}$ is the unique point on $\gamma_{1}$ which has the same
 fate as $x$:that is their trajectory are asymptotic with the rate
 of $exp(-t)$, see [1],Chapter $13$\\
 The  Integrals  converge since the  corresponding functions
 approach to zero with $"exp"$ rate,
 For x on $\gamma_{1}$ ,$x^{*}=x$ and g was an integral for f
 restricted on $\gamma_{1}$ \\
 This Shows that g described above is an integral for f in a
 neighborhood of limit cycle  $\gamma_{1}$ with the same values on
 $\gamma_{1}$.

 g can be define  on the whole of the  plane since the orbits of exterior points of $\gamma_{n}$
 accumulate to it.\\
 Note  that Since  Vector  field  $L$ is analytic,$x^{*}$ Is  analytic too,thus
 the  proposition is valid if we define  L on smooth or analytic
 function space\\

 \textbf{Remark 1.}Let  F be an even polynomial,then the corresponding
 Lienard equation L, has a center and both kernel and co-kernel
 of operator $L(f)=L.f$ are infinite dimensional space:\\
 We Show that L possesses a global analytic  first integral thus
 kernel of above operator is infinite  dimensional space,furthermore  for each set of n closed orbits we present $n$ independent
 elements in the  quotion space of Image of operator $L(f)=L.f$:
 let f be a smooth (analytic or algebraic) maps separates
 closed the  orbits then the elements $1,f,f^{2},f^{3}.....f^{n-1}$ are
 independent in quotion space of image because for each map g
 ,$L.g$ should vanish in  at least one point on  each closed
 orbit\\
 It remains to prove the integrability of the Lienard equation with
 center:Let $F(x)=K(x^{2})$ for a polynomial K,the square of
 intersection of orbits with the graph of parabola $y=-x^{2}$
 defines a global first integral,this parabola is not transverse  to
 lienard vector field  at the  origin so apparently the above first integral is
 not analytic at the origin.Using Change of coordinate $x:=x^{2}$
 $y:=y$ we find that this first integral corresponds  to intersection
 of solutions of
  $$
\begin{cases}
 \dot{x}=y-k(x)\\
 \dot{y}=-1
\end{cases}
$$
with transversal section $y=x$,which is analytic\\

\textbf{ Remark 2.}The operator L described above can be
restricted to algebraic functions.Since for each set of n closed
orbits we can give n independent element in co-kernel,The index of
operator is an upper bound for the number of the limit cycles.In
line of conjecture in $[2]$ on the number of limit cycles of
lienard equation,we suggest to compute the index of operator
restricted to polynomials maps.Are there uniform upper bounds for
this index in
terms of degree of F,where the degree of  F is odd\\

\textbf{ Remark 3.}The $Fredholm$ Index ,mentioned above, is  not
necessarily finite if an arbitrary algebraic Vector Field Possess
A limit cycle,for example there is a cubic system with a center
and a limit cycle simultaneously ,so in this case the co-kernel's dimension is infinite.\\
But not only such "co-existence" of limit cycle and center cannot
occur   for quadratic Systems,but also,all  quadratic systems with
center have been classified with a FINITE number of algebraic
condition.Furthermore,since $fredholm$ index is fixed on the
connected component of the space of $fredholm$ operators ,it
strongly  seems that  $H(2)$ is finite provided that  we prove
following
problem :\\

\textbf{Problem.}For a quadratic system without center the
corresponding operator $X(f)=X.f$ is $fredholm$ where X is
Considered as an analytic Vector field on 2 dimensional sphere(X
is considered as poincare compctification of corresponding
quadratic system )

\end{document}